\documentclass[a4paper, 11pt]{article}
\usepackage{latexsym} 
\usepackage{amsfonts} 
\usepackage{amsmath} 
\usepackage{amsthm} 
\usepackage{mathrsfs} 
\usepackage{dsfont} 
\usepackage{bbold} 
\usepackage[english]{babel} 
\usepackage{caption} 
\usepackage{epsfig} 
\usepackage{float} 
\usepackage{subfigure} 
\usepackage{psfrag} 
\usepackage{graphicx} 
\usepackage{epsfig} 
\usepackage{amsbsy} 
\usepackage{hyperref}
 
\DeclareMathOperator{\Val}{\matV}

\newtheorem{theorem}{Theorem} 
\newtheorem*{prop*}{Theorem} 
\newtheorem*{hyp*}{Hypothesis} 
\newtheorem{thm}{Theorem}

\newtheorem{lemma}[theorem]{Lemma}

\newtheorem{hyp}{Hypothesis}

\newcommand{\zerarcounters}{\setcounter{equation}{0}\setcounter{theorem}{0}}

\newcommand{\ZZZ}{\mathds{Z}} 
\newcommand{\CCC}{\mathds{C}} 
\newcommand{\NNN}{\mathds{N}} 
 
\newcommand{\RRR}{\mathds{R}} 
\newcommand{\TTT}{\mathds{T}} 
\newcommand{\uno}{\mathds{1}}

\newcommand{\CCCC}{{\mathcal C}}

\newcommand{\calG}{{\mathcal G}} 
\newcommand{\calH}{{\mathcal H}}

\newcommand{\TT}{{\mathcal T}}

\newcommand{\gotC}{{\mathfrak C}}

\newcommand{\matA}{{\mathscr A}}

\newcommand{\matV}{{\mathscr V}}

\newcommand{\Fullbox}{{\rule{2.0mm}{2.0mm}}} 
 
\newcommand{\EP}{\hfill\Fullbox\vspace{0.2cm}} 
\newcommand{\prova}{\noindent{\it Proof. }} 
\newcommand{\io}{\infty} 
\newcommand{\e}{\varepsilon} 
\newcommand{\al}{\alpha} 
\newcommand{\de}{\delta} 
\newcommand{\be}{\beta}

\newcommand{\x}{\xi} 
 
\newcommand{\ka}{\kappa} 
 
\newcommand{\om}{\omega}

\newcommand{\oo}{\omega}

\newcommand{\nn}{\nu} 
\newcommand{\pps}{\psi} 
\newcommand{\vzero}{0} 
\newcommand{\xx}{x}

\newcommand{\cc}{c}

\newcommand{\ii}{{\rm i}}

\def\tilde#1{\widetilde{#1}}
 
\def\ins#1#2#3{\vbox to0pt{\kern-#2 \hbox{\kern#1 #3}\vss}\nointerlineskip}

\begin{document}

\title{\bf Forced quasi-periodic oscillations in strongly dissipative systems of any finite dimension}
\author 
{\bf Guido Gentile, Alessandro Mazzoccoli, Faenia Vaia
\vspace{2mm} 
\\ \small
Dipartimento di Matematica, Universit\`a di Roma Tre, Roma, I-00146, Italy
\\ \small 
E-mail: gentile@mat.uniroma3.it, alessandro.mazzoccoli@gmail.com, faenia.vaia@uniroma3.it
}

\date{} 
 
\maketitle 
 
\begin{abstract} 
We consider a class of singular ordinary differential equations
describing analytic systems of arbitrary finite dimension,
subject to a quasi-periodic forcing term and in the presence of dissipation. 
We study the existence of response solutions,
i.e. quasi-periodic solutions with the same frequency vector as the forcing term,
in the case of large dissipation. We assume the system to be conservative in the absence
dissipation, so that the forcing term is --- up to the sign --- the gradient of a potential energy,
and both the mass and damping matrices to be symmetric and positive definite.
Further, we assume a non-degeneracy condition on the forcing term, essentially that
the time-average of the potential energy has a strict local minimum.
On the contrary, no condition is assumed on the forcing frequency;
in particular we do not require any Diophantine condition.
We prove that, under the assumptions above, a response solution always exist provided
the dissipation is strong enough.
This extends results previously available in literature in the one-dimensional case.
\end{abstract} 

\zerarcounters 
\section{Introduction} 
\label{sec:1} 

Consider the singular non-autonomous ordinary differential equation on $\RRR^m$
\begin{equation} \label{eq:1.1}
\e \ddot x + \Gamma \dot x + \e g(x) = \e f(\oo t) ,
\end{equation}
where $x=(x_1,\ldots,x_m)$, the dots denote derivatives with respect to time $t$,
$\oo \in \RRR^{d}$ is the frequency vector of the forcing term, $\e \in \RRR_{+}$
is a small parameter and $\Gamma$ is a positive definite symmetric $m \times m$ real matrix.
The functions $f: \TTT^{d} \to \RRR^{m}$ and $g: \RRR^{m} \to \RRR^{m}$
are assumed to be real analytic, with $f$ quasi-periodic in $t$,
so we can write its Fourier series expansion as
\begin{equation} \nonumber
f(\pps)= \sum_{\nn \in \ZZZ^{d}} e^{\ii \nn \cdot \pps} f_{\nn} ,
\end{equation}
where $\pps \in \TTT^{d}$ and $f_{\nn}$ are the Fourier coefficients.  Here and henceforth $\cdot$ denotes
the standard scalar product in $\RRR^d$, i.e. $\nn\cdot\pps=\nu_1\psi_1+\ldots+\nu_d\psi_d$.
Since the results we are going to discuss are local, we could assume $g$ to be defined only in an open set
$\matA\subset\RRR^{m}$.

Systems described by equations of the form \eqref{eq:1.1} naturally arise in classical mechanics and
electronic engineering \cite{S,AA,A,Gans,ADH,DM,BDGM}. As a mechanical system, \eqref{eq:1.1} describes
a collection of points of unit mass; as we show in Section \ref{sec:6}, the analysis easily extends to the case
of points with different masses. The points interact through a mechanical force $-g$ and are subject to a forcing term $f$,
while the term proportional to the velocity (drag) takes into account dissipation effects --- described
as dashpots or dampers --- due, for instance, to air resistance, fluid viscosity or friction in the gears of a device. 
Thus, $\Gamma$ is called the \emph{damping matrix}. In \eqref{eq:1.1}
we assume the matrix $\Gamma$ to be a positive definite symmetric matrix, as such a property is satisfied in many problems
of physical relevance arising in fluid dynamics, electrodynamics and mechanics \cite{A,D}.
While in principle one can also consider situations in which the matrix $\Gamma$ is not symmetric \cite{AA},
the results we shall present in this paper hold for symmetric damping matrices.
The parameter $\e$ plays the role of a \emph{perturbation parameter} and measures the intensity of the dissipation
(for this reason we take it to be positive): small $\e$ means strong dissipation.

We also assume that the system is conservative in the absence of dissipation,
so that there exists a real analytic function $V\!:\RRR^{m} \to \RRR$ such that
\begin{equation}
\label{eq:1.2}
g(x) = \frac{\partial V}{\partial x}(x) ,
\end{equation}
i.e. the mechanical force $-g$ is a gradient and $V$ is the potential energy.

We are interested in studying the existence, for small values of $\varepsilon$, of \emph{response solutions} to \eqref{eq:1.1},
i.e. quasi-periodic solutions with the same frequency vector $\oo$ as the forcing \cite{CLB}.
If $\varepsilon=0$, for any constant vector $\cc\in \RRR^{m}$, $\xx=\cc$ is a solution to \eqref{eq:1.1}.
The problem we consider is whether it is possible to choose the constant vector $\cc$ so that,
for $\e$ small enough, the equation \eqref{eq:1.1} admits a response solution which tends to $\cc$ as $\e$ tends to zero.

The case $d=1$ has been already studied in the literature. 
Without any assumptions on the functions $f$ and $g$, the answer in general is negative \cite{G1}.
The analysis in \cite{GV} shows that, 
for any frequency vector $\oo$, if $\cc$ is a simple zero of the the function $h(x):=g(x)-f_{\vzero}$,
then for $\e$ small enough there is a response solution which tends to $\cc$ as $\e$ tends to zero.
This strengthens previous results, where non-resonance conditions were assumed on $\oo$ \cite{GBD1,CCD,CFG1,CFG2}. 
More or less strong non-resonance conditions assure some regularity in the dependence on $\e$ of the solution;
in particular, without assuming any condition on $\oo$, in general no more than a continuous dependence on $\e$
can be obtained, while smoothness and analyticity -- and even Borel-summability -- in domains tangential to the imaginary axis
follow if suitable Diophantine conditions are assumed on $\oo$ \cite{GBD2,CCD,CFG1}.
On the other hand, on physical grounds, in the presence of dissipation one may expect a response solution to exist
independently on the arithmetic properties of the frequency vector. However, from a mathematical point of view,
quasi-periodicity typically involves small divisor problems, and to prove the existence of quasi-periodic
solutions for any rotation vector may be a rather delicate problem, also in the context of singularly perturbed systems
such as \eqref{eq:1.1}; see \cite{GV} for further references and comments.

Existence of quasi-periodic solutions in singular ordinary differential equations
has been extensively studied in the literature. Usually one makes some non-degeneracy assumption
on the unperturbed solution, which ensures that the response solution either smoothly appears
by bifurcation, due to the loss of stability of an isolated solution, such as an equilibrium point or a periodic solution,
or persists as a stable or hyperbolic central manifold; see for instance \cite{A,BC,BB,SPR}.
By contrast, in \eqref{eq:1.1} any constant is a solution at $\e=0$, 
so that the unperturbed bifurcating solution is not given a priori; in particular it does not satisfy
any stability or hyperbolicity condition.

In this paper, we extend the result in \cite{GV} to the higher dimensional case $m>1$.
Thus, we may and shall only require that the components of $\oo$ are rationally independent, i.e.
$ \oo \cdot \nn \ne 0 \, \, \forall \nn \in \ZZZ_{*}^{d}:= \ZZZ^{d} \setminus {\{\vzero\}}$;
this is not restrictive, because one can always reduce to such a case, possibly changing the value of $m$.
Like in the one-dimensional case,
some non-degeneracy condition involving the functions $f$ and $g$ must be required. Let us define the function
\begin{equation}
\label{eq:1.3}
U(x) := V(x) - \langle f_0 , x \rangle ,
\end{equation}
where $V$ is the potential energy associated to $g$ through \eqref{eq:1.2}
and $\langle \cdot , \cdot \rangle$ denotes the standard scalar product in $\RRR^m$.

\begin{hyp}
\label{hyp:1}
The function $U(x)$ in \eqref{eq:1.3} has a strict local minimum point $c$.
\end{hyp}

We shall prove the following result.

\begin{thm} 
\label{thm:1}
Consider the ordinary differential equation \eqref{eq:1.1}, with $g$ as in \eqref{eq:1.2},
and assume Hypothesis \ref{hyp:1}.
For any frequency vector $\oo \in \RRR^{d}$, there exists $\e_{0}>0$ such that for all $\e \in (0,\e_{0})$
there is a quasi-periodic solution $x_{0}(t)=c + X(\oo t,\e)$ to \eqref{eq:1.1},
such that $X(\pps,\e)$ is analytic in $\pps$ and goes to $0$ as $\varepsilon \to 0^{+}$.
\end{thm}

Under the assumption that $\oo$ satisfies the Bryuno condition (a weaker non-resonance condition, 
with respect to the standard Diophantine condition), the existence of response solutions
has been proved in \cite{M},
by extending the results in \cite{G2} to higher dimensions. In this paper, we follow
the approach introduced in \cite{GV} to remove the non-resonance condition on the frequency vector.

What we really need about the function $U$ is that its Hessian matrix is strictly definite.
Thus, the same result of existence holds if, instead of Hypothesis \ref{hyp:1}, we assume $\cc$ to be
a strict local maximum for $\cc$. 
Finally, we might allow $\e$ to be negative as well, but we prefer to fix $\e>0$ because this
corresponds to a physical system in the presence of dissipation.

The results above can be extended to differential equations equations of the form
\begin{equation} 
\label{eq:1.4}
\e M \ddot x + \Gamma \dot x + \e h(x,\oo t) = 0 , 
\end{equation}
where $M$, the mass matrix, is a diagonal or, more generally, a positive definite symmetric matrix, and
$h:\RRR^{m}\times \TTT^d \to \RRR^m$ is a real analytic function. For instance,
if $m=3n$ and $M$ is diagonal, with $M_{ii}=M_1$ for $i=1,2,3$, $M_{ii}=M_2$ for $i=4,5,6$,
\ldots, $M_{ii}=M_n$ for $i=3n-2,3n-1,3n$, then \eqref{eq:1.4} describes a system of $n$ points
in $\RRR^{3}$, of masses $M_1,M_2,\ldots,M_n$, respectively.
We assume $h$ to be the gradient of a time-dependent
potential energy $\tilde{V}$, i.e.
\begin{equation}
\label{eq:1.5}
h(x,\oo t) = \frac{\partial \tilde{V}}{\partial x}(x,\oo t) .
\end{equation}
Let us consider the Fourier expansion of the function $\tilde{V}$:
\begin{equation}
\label{eq:1.6}
\tilde{V}(x,\pps)= \sum_{\nn \in \ZZZ^{d}} e^{\ii \nn \cdot \pps} \tilde{V}_{\nn}(x) .
\end{equation}
Then we assume the following condition on $\tilde{V}$.

\begin{hyp}
\label{hyp:2}
The function $\tilde{V}_0(x)$ in \eqref{eq:1.6} has a strict local minimum point $c$.
\end{hyp}

Then the following result generalises Theorem \ref{thm:1}.

\begin{thm} 
\label{thm:2}
Consider the ordinary differential equation \eqref{eq:1.4}, with $h$ as in \eqref{eq:1.5},
assume Hypothesis \ref{hyp:2}.
For any frequency vector $\oo \in \RRR^{d}$, there exists $\e_{0}>0$ such that for all $\e\in(0,\e_{0})$
there is a quasi-periodic solution $x_{0}(t)=c + X(\oo t,\e)$ to \eqref{eq:1.4},
such that $X(\pps,\e)$ is analytic in $\pps$ and goes to $0$ as $\varepsilon \to 0^{+}$.
\end{thm}

The rest of the paper is organised as follows.
We first discuss the proof of Theorem \ref{thm:1}. In Section \ref{sec:2} we outline the overall strategy
of the proof, by showing that the problem can be split into two different equations: the range equation
and the bifurcation equation. In Section \ref{sec:3} we prove the existence of a solution to the
range equation depending on a parameter $\zeta$. In Section \ref{sec:4} we fix the parameter
$\zeta$ as a function of $\e$ in such a way that the bifurcation equation is solved as well.
This completes the proof of Theorem \ref{thm:1}.
Finally, in Section \ref{sec:6} we discuss how to generalise the results and to prove Theorem \ref{thm:2}.

On physical grounds, one may expect the response solution to be attracting. As a byproduct,
local uniqueness would also follows. In the one-dimensional case, this can be easily
proved by a rescaling of time \cite{BDG,GV}. However, the argument does not extend
to higher dimensions. So, uniqueness the response solutions
described in Theorems \ref{thm:1} and \ref{thm:2} is left as an open problem.

Another issue which deserves further investigation is how to deal with systems in which
the mass matrix and the damping matrix are not symmetric.
From a physical point of view, it is not restrictive to assume the mass matrix $M$ to be symmetric
since the kinetic energy is a quadratic form; more generally $M$ may depend on the coordinates $x$,
but the analysis below can be extended to cover such a case.
As far as to the damping matrix, taking $\Gamma$ to be positive means that we are assuming
that the friction that a mass point experiences because of the presence of another mass point
is the same as that it produces on that point. Usually, this is what happens in most dynamical structures \cite{D};
however, there are situations in which the damping matrix are asymmetric \cite{AA}.
From a technical point of view, the property of $\Gamma$ being symmetric is used to control
the inverse of a suitable matrix which appears along the analysis: renouncing such a property
would introduce severe technical difficulties.
This is not surprising, since such kind of problems are known to arise
in the higher dimensional case -- see for instance \cite{EK,GP,CG} and references therein.
Similar considerations apply to the case in which the force is not conservative.

\zerarcounters 
\section{Plan of the proof of Theorem \ref{thm:1}}
\label{sec:2}

We first prove Theorem \ref{thm:1}: here we describe the strategy of the proof, with
the details to be worked out in the forthcoming sections \ref{sec:3} and \ref{sec:4}.
Eventually, in Section \ref{sec:6}, we discuss
how to adapt the proof so as to prove the more general Theorem \ref{thm:2}.

Let us denote by $\Sigma_{\xi}$ the strip around $\TTT^{d}$ of width $\xi$ and
by $B_{\rho}(c)$ the disk of center $c=(c_1,\ldots,c_m)$ and radius $\rho$ in the $m$-dimensional complex plane,
\begin{equation} \nonumber
\Sigma_{\xi} := \left\{ \pps \in \CCC^d : {\rm Re}\,\psi_i \in \TTT, \quad |{\rm Im}\psi_i| < \xi \right\} ,
\qquad B_{\rho}(c) := \left\{  x \in \CCC^m : |x-c| < \rho \right\} .
\end{equation}
By the assumptions on $f$ and $g$, for any $c \in\RRR^{m}$ there exist
$\xi_{0}>0$ and $\rho_{0}>0$ such that $\pps \mapsto f(\pps)$ is analytic
in $\Sigma_{\xi_{0}}$ and $x \mapsto g(x)$ is analytic in $B_{\rho_{0}}(c)$.
Then for all $\xi<\xi_{0}$ one has
\begin{equation} \label{eq:2.1}
f_i(\pps) = \sum_{\nn\in\ZZZ^{d}} {\rm e}^{\ii\nn\cdot\pps} (f_{\nn})_i , \qquad 
|(f_{\nn})_i| \le \Phi \, {\rm e}^{-\xi |\nn|} , \qquad i=1,\ldots,m,
\end{equation}
where $\Phi$ is the maximum of $|f(\pps)|$ for $\pps\in\Sigma_{\xi}$.

Define $G_p(c)$ as the tensor of rank $p+1$ with components 
\begin{equation} \nonumber 
[G_{p}(c)]_{i,i_1,\ldots,i_p} = \left. \frac{1}{p!} \frac{\partial^{p}}{\partial x_{i_1} \ldots \partial x_{i_p} } g_i (x) 
\right|_{x=c} ,
\end{equation}
where $i,i_1,\ldots,i_p=1,\ldots,m$, and write
\begin{equation} \nonumber
g_i(x)= \sum_{p=0}^{\io} [G_{p}(c) (x-c)^{p}]_{i} = \sum_{p=0}^{\io}
\sum_{i_{1},\ldots,i_{p}=1}^{m} [G_{p}(c)]_{i,i_1,\ldots,i_p} \prod_{j=1}^{p} \left( x_{i_j} - c_{i_j} \right) , 
\qquad i=1,\ldots,m .
\end{equation}
For all $\rho < \rho_{0}$ one has
\begin{equation} \label{eq:2.2}
\left| [G_{p}(c)]_{i,i_1,\ldots,i_p} \right| \le \Delta \rho^{-p} ,
\end{equation}
with $\Delta$ being the maximum of $|g(x)|$ for $x \in B_{\rho_0-\rho}(c)$. 

Let us rewrite (\ref{eq:1.1}) as
\begin{equation} \label{eq:2.3}
\varepsilon \ddot x + \Gamma \dot x + \e \, g(c) + \e \, A \left( x-c \right) + \varepsilon \, G(x) = \varepsilon \, f(\oo t) ,
\end{equation}
where $A=G_1(c)$ is the $m \times m$ matrix with entries
\begin{equation} \nonumber
A_{i,j} := \left. \frac{\partial g_i}{\partial x_j} (x) \right|_{\xx=\cc} =
\left. \frac{\partial^2 V}{\partial x_i \partial x_j} (x) \right|_{\xx=\cc} ,
\qquad i,j =1,\ldots, m ,
\end{equation}
and
\begin{equation} \nonumber
G(x) := g(x)-g(c)- A \left( x - c \right) = \sum_{p=2}^{\io} G_{p}(c) (x-c)^{p} ,
\end{equation}
By Hypothesis \ref{hyp:1}, the symmetric matrix $A$ is positive definite; hence its eigenvalues are all positive.

We look for a quasi-periodic solution to \eqref{eq:2.3}, that is a solution of the form
\begin{equation} \label{eq:2.4}
x(t,\e) = c + \zeta + u(\oo t, \e,\zeta) , \qquad u(\pps,\e,\zeta) = 
\sum_{\nn\in\ZZZ^{d}_{*}} {\rm e}^{\ii\nn\cdot\pps} u_{\nn} ,
\end{equation}
where $\zeta$ is a real parameter that has to be fixed
eventually, and $\pps \mapsto u(\pps,\varepsilon,\zeta)$
is a zero-average quasi-periodic function, with Fourier coefficients depending
on both $\varepsilon$ and $\zeta$.

It is more convenient to write \eqref{eq:2.3} in Fourier space, so by choosing $c$ as in Hypothesis \ref{hyp:1}, 
%
%
we obtain the following equations:
\begin{subequations} \label{eq:2.6}
\begin{align}
\left( - \e \, (\oo\cdot \nu)^{2} \uno + \ii \oo \cdot \nu \Gamma + \e A \right)
u_{\nn} & = - \varepsilon \left[ G(c + \zeta+u) \right]_{\nn} + \varepsilon \, f_{\nn} , \qquad \nn \neq \vzero ,
\label{eq:2.6a} \\
\varepsilon \, A \, \zeta  & = - \varepsilon \left[ G(c + \zeta+u) \right]_{\vzero} ,
\label{eq:2.6b}
\end{align}
\end{subequations}
where $\uno$ is the $m \times m$ identity matrix and the notation $[G(c+\zeta+u)]_{\nn}$ means that we first write $u$ according to \eqref{eq:2.4},
then expand $G(c+\zeta+u)$ in Fourier series in $\pps$ and finally keep the Fourier coefficient
with index $\nn$. We call \eqref{eq:2.6a} the \emph{range equation} and \eqref{eq:2.6b} 
the \emph{bifurcation equation} \cite{CH,GS,SLSF}.

To simplify the notation, we call
\begin{equation} \label{eq:2.5}
D(\e,s) := - \e \, s^{2} \uno + \ii s \Gamma + \e A
\end{equation}
and write \eqref{eq:2.6a} as
\begin{equation} \nonumber
D(\e,\oo\cdot\nu) \, u_{\nn} = - \varepsilon \left[ G(c + \zeta+u) \right]_{\nn} + \varepsilon \, f_{\nn} , \qquad \nn \neq \vzero .
\end{equation}

In order to study the equations \eqref{eq:2.6}, we start by ignoring the bifurcation equation \eqref{eq:2.6b} and
considering only the range equation \eqref{eq:2.6a}. We treat $\zeta$ as a free parameter,
close enough to $0$, and look for a solution depending on such a parameter, by using perturbation theory.
To this aim, we modify \eqref{eq:2.6a} introducing the auxiliary parameter $\mu$,
\begin{equation} \label{eq:2.7}
D(\varepsilon,\oo\cdot\nn) \, u_{\nn} = - \mu \left[ \varepsilon G(c+\mu\zeta + u) \right]_{\nn} + 
\mu \varepsilon \, f_{\nn} , \qquad \nn \neq \vzero .
\end{equation}
Of course the original equation \eqref{eq:2.6a} is recovered only for $\mu=1$, so eventually we have to check
that the radius of convergence is larger than 1.

Then, in Section \ref{sec:3}, we show that a quasi-periodic solution to \eqref{eq:2.7} exists in the form of a power series in $\mu$,
\begin{equation}  \label{eq:2.8}
u(\oo t, \varepsilon,\zeta,\mu) = \sum_{k=1}^{\infty} \sum_{\nn\in\ZZZ^{d}_{*}} 
\mu^{k} {\rm e}^{\ii\nn\cdot\pps} u^{(k)}_{\nn} ,
\end{equation}
and we prove that the radius of convergence of \eqref{eq:2.8} is strictly greater than one.

Once a solution $u(\om t , \e, \zeta) := u(\om t,\e,\zeta,1)$ to the range equation has been proved to exist,
we pass to the bifurcation equation \eqref{eq:2.6b}. The latter is studied in Section \ref{sec:4}
as an implicit function problem to determine the constant $\zeta$ in terms of $\e$. The solution provides 
a value $\zeta=\zeta(\e)$, depending continuously on $\e$, 
such that $x(t,\e)=c + \zeta(\e)+u(\om t,\e,\zeta(\e),1)$ is a solution to \eqref{eq:2.3}
and hence to \eqref{eq:1.1}.

This completes the proof of Theorem \ref{thm:1}.

\zerarcounters 
\section{The range equation}
\label{sec:3}

In this and in the following two sections we provide the details of the proof of Theorem \ref{thm:1}.
Here we show that, for $\e$ and $\zeta$ small enough, there exists a solution of the
form \eqref{eq:2.4} to the range equation \eqref{eq:2.6a}. Moreover, such a solution
$u(\oo t,\e,\zeta)$ is analytic in $\pps=\oo t$ and $\zeta$, and tends to $c$ as $\e$ tends to $0$.

\subsection{Recursive equations}

The coefficients $u^{(k)}_{\nn}$ are obtained recursively
by inserting the expansion \eqref{eq:2.8} into \eqref{eq:2.7}. If we set $u^{(1)}_{\vzero}=\zeta$
and $u^{(k)}_{\vzero}=\vzero$ $\forall k\ge 2$, to make the notations uniform,
then, one obtains, formally,
\begin{subequations}  \label{eq:2.9}
\begin{align}
D(\varepsilon, \oo\cdot\nn) \, u^{(1)}_{\nn} & = \varepsilon \, f_{\nn} 
\label{eq:2.9a} \\
D(\varepsilon, \oo\cdot\nn) \, u^{(k)}_{\nn} & = - \varepsilon \, 
\sum_{p=2}^{\io} G_{p}(c) \!\!\!\!\!\!
\sum_{\substack{ k_{1},\ldots,k_{p} \ge 1 \\ k_{1}+\ldots+k_{p}=k-1}}
\sum_{\substack{\nn_{1},\ldots,\nn_{p} \in\ZZZ^{d} \\ \nn_{1}+\ldots+\nn_{p}=\nn}} \!\!\!\!\!\!
u^{(k_{1})}_{\nn_{1}} \ldots u^{(k_{p})}_{\nn_{p}} , \qquad k \ge 2 ,
\label{eq:2.9b}
\end{align}
\end{subequations}
where $\nn\neq \vzero$. We assume that the sums over the empty set are meant as zero.
It is easily checked that $u^{(2)}_{\nn}=0$ $\forall\nn\in\ZZZ^{d}_{*}$.

\subsection{Tree representation for the coefficients}

By iterating \eqref{eq:2.9b}, one obtains a diagrammatic representation of the coefficients
$u^{(k)}_{\nn}$ in terms of trees; we refer to \cite{G3,GV} and references therein for further details.

A \textit{rooted tree} $\theta$ is a graph with no cycle,
such that all the lines are oriented toward a unique
point, called the \textit{root}, which has only one incident line, called the root line.
All the points in $\theta$ except the root are called \textit{nodes}.
Given a rooted tree $\theta$ we denote by $N(\theta)$ the set of its nodes
and by $L(\theta)$ the set of its lines.

The orientation of the lines in $\theta$ induces a partial ordering 
relation $\preceq$ between the nodes. If $v$ and $w$ are two nodes of the tree,
we write $w \prec v$ every time $v$ is along the path
of lines which connects $w$ to the root; we shall write $w\prec \ell$ if
$w\preceq v$, where  $v$ is the unique node that the line $\ell$ exits.
For any node $v$ let $p_{v}$ denote the number of lines entering $v$.
One has $p_v \ge 0$; if $p_v=0$ we say that $v$ is an \emph{end node},
while if $p_v\ge 1$ we say that $v$ is an \emph{internal node}.
Denote by $E(\theta)$ the set of end nodes and by $V(\theta)$
the set of internal nodes in $\theta$; one has $N(\theta)=E(\theta) \amalg V(\theta)$.

For any discrete set $A$, we denote by $|A|$ its cardinality;
we define the \textit{order} of $\theta$ as $k(\theta):=|N(\theta)|$.

We associate with each end node $v\in E(\theta)$ a \textit{mode} label $\nn_{v}\in\ZZZ^{d}$.
According the the values of the mode labels, we split $E(\theta)$ into two disjoint sets by writing
$E(\theta)=E_{0}(\theta) \amalg E_{1}(\theta)$, with
$E_{0}(\theta)=\{ v \in E(\theta) : \nn_{v} = \vzero \}$ and
$E_{1}(\theta)=\{ v \in E(\theta) : \nn_{v} \neq \vzero \}$.
With each line $\ell\in L(\theta)$ we associate
a \textit{momentum} $\nn_{\ell} \in \ZZZ^{d}$ with the constraint
\begin{equation} \label{eq:3.1}
\nn_{\ell}=\sum_{\substack{w\in E(\theta) \\ w \prec \ell}} \nn_{w} .
\end{equation}
Finally we impose the further constraints that
\begin{itemize}
\itemsep0em
\item $p_{v}\ge 2$ $\forall v\in V(\theta)$,
\item $\nn_{\ell} \neq \vzero$ for any line $\ell$ exiting a node in $V(\theta)$.
\end{itemize}

We call \textit{equivalent} two labelled rooted trees which can be transformed into
each other by continuously deforming the lines in such a way that they do not cross each other. 
Since, in the following, we shall consider only inequivalent labelled rooted trees, 
for simplicity's sake they will be called simply `trees'.

Next, with each node $v\in N(\theta)$  we associate a \textit{node factor}
\begin{equation} \label{eq:2.10}
F_{v} := \begin{cases}
- \e \, G_{p_{v}} , & v \in V(\theta) , \\
\e \, f_{\nn_{v}} , & v \in E_{1}(\theta) , \\
\zeta , & v \in E_{0}(\theta) ,
\end{cases}
\end{equation}
and with each line $\ell\in L(\theta)$ we associate a \textit{propagator}
\begin{equation} \label{eq:2.11}
\mathcal{G}_{\ell} := \begin{cases}
D^{-1}(\e, \oo\cdot\nn_{\ell}) , & \nn_{\ell} \neq \vzero , \\
\uno , & \nn_{\ell}=\vzero .
\end{cases}
\end{equation}
Note that the node factor \eqref{eq:2.10} is a tensor of rank $p_{v}+1$, while the propagator \eqref{eq:2.11} is a matrix.
Finally we define the value of the tree $\theta$ as
\begin{equation} 
\label{eq:3.2}
\Val(\theta) := \Biggl( \prod_{v\in N(\theta)} F_{v} \Biggr) 
\Biggl( \prod_{\ell\in L(\theta)} \mathcal{G}_{\ell} \Biggr) ,
\end{equation}
where
\begin{itemize}
\itemsep0em
\item the indices of both node factors and propagators are such that, if $[F_v]_{i_0,i_1,\ldots,i_{p_v}}$
is the component of $F_v$, then the second index of the propagator of the line $\ell_v$ is $i_0$ 
and the first indices of the propagators of the $p_v$ lines entering $v$ are $i_1,\ldots,i_{p_v}$;
\item all the indices have been contracted, so that \eqref{eq:3.2} gives a vector,
whose component index $i$ equals the first index of the propagator of the root line.
\end{itemize}

If we denote by $\TT_{k,\nn}$ the set of non-equivalent trees of order $k$
and momentum $\nn$ associated with the root line, and define the coeffiecents $u^{(k)}_{\nu}$ in \eqref{eq:2.8} as
\begin{equation} 
\label{eq:3.3}
u^{(k)}_{\nn} := \sum_{\theta\in \TT_{k,\nn}} \Val(\theta) , \qquad \nn \in\ZZZ^{d}_{*} ,
\end{equation}
then it is straightforward to see that the equations \eqref{eq:2.9} are solved for all $k\in\NNN$.

\subsection{Bounds on the propagators}
\label{sec:3.03}

We can estimate the norm of the propagators as follows.
By construction $\Gamma$ is symmetric
and positive definite, so we write $\Gamma:= K^{2}$, with $K$ a symmetric positive definite matrix.
Let us denote by $\ka_{1} \le \dots \le \ka_{m}$ the eigenvalues of $K$.
Hence, we can write $D=D(\varepsilon, s)$ as
\begin{equation} \nonumber
D = K \left( \ii s \uno + K^{-1} (-\e s^2 \uno + \e A) K^{-1} \right) K ,
\end{equation}
so that
\begin{equation} \nonumber
D^{-1}=K^{-1 } \left( \ii s \uno + K^{-1} (-\e s^2 \uno + \e A) K^{-1} \right)^{-1} K^{-1}.
\end{equation}
By defining 
\begin{equation} \nonumber
S := \e \, K^{-1} A K^{-1} - \e s^{2} K^{-2} \quad \text{and} \quad R:= S + \ii s \uno ,
\end{equation}
we obtain
%
$D^{-1}= K^{-1} R^{-1} K^{-1}$.
%
Let $b_{1} \le \dots \le b_{m}$ be the eigenvalues of $K^{-1}AK^{-1}$. Since both $A$ and $K$ are positive definite,
$b_1,\ldots,b_m$ are positive. By Rellich's theorem \cite{R},
it follows that the eigenvalues of $S$ are $b_{k}\e + O(\e s^{2})$, $k=1,\dots, m$, and hence
the eigenvalues of $R$ are
$\be_{k} + \ii s$, $k=1,\dots, m$, with $\be_{k} = b_{k}\e + O(\e s^{2})$.
%
%
Set $\be:=\min\{|\be_1|,\ldots,|\be_m|\}$.

Let us denote by $|\cdot|$ the Euclidean norm on $\RRR^m$ and by $\| R \| := \sup_{|x|=1} |Rx|$ the induced norm;
one has
\begin{equation} \nonumber
\|R^{-1}\| \le \frac{1}{\sqrt{\be^{2} + |s|^{2}}}.
\end{equation}
For any $s \in \RRR$ one has $|\be_{k} - b_{k}\e | \le W |\e s^{2}|$ for a suitable constant $W$.
Thus, if $|s| \le \al$, where $\al := \sqrt{b_{1}/2W}$, one obtains
\begin{equation} \nonumber
\|R^{-1}\| \le \min \left\{\frac{2}{|b_{1}\e|},\frac{1}{|s|}\right\} ,
\end{equation}
otherwise one has $\|R^{-1}\| \le 1/\al$. 
By collecting together the bounds above, we arrive at
\begin{equation} \label{eq:3.5a}
\| D^{-1} \| \le \frac{1}{\ka^{2}_{1}} \max \left\{ \frac{1}{\alpha} , \min
\left\{ \frac{1}{|s|}, \frac{2}{|b_{1}\e|} \right\} \right\} .
\end{equation}
Since we are interested in results for $\e$ small enough, in the following we assume $\e<\e_1$, 
where $\e_1:=\al/b_1$, so that we can replace \eqref{eq:3.5a} with
\begin{equation} \label{eq:3.5}
\| D^{-1} \| \le \frac{1}{\ka^{2}_{1}}  \min \left\{  \frac{2}{|b_{1}\e|} , \max  \left\{ \frac{1}{\al} , \frac{1}{|s|} \right\} \right\} .
\end{equation}
%

\subsection{Bounds on the coefficients}
\label{sec:3.4}

For future convenience, set
\begin{equation} \label{eq:3.6}
C_0 := m^{2}\rho^{-1}\max \left\{ \Phi, 1, \frac{2\Delta}{\ka_{1}^{2}|b_{1}|}\right\} ,
\end{equation}
and define
\begin{equation} \label{eq:3.6a}
s_{N} : = \min\{ |\oo\cdot\nn| : 0<|\nn| \le N \} , \qquad r_{N} = \min \{ s_{N} , \al \} .
\end{equation}

By reasoning as in the proof of Lemma 2.3 in \cite{CFG2}, for any tree $\theta$
one obtains the bound 
\begin{equation} \label{eq:3.4}
|E(\theta)|\ge \frac{1}{2} \left( k(\theta)+1 \right) .
\end{equation}
%

\begin{lemma} 
\label{lem:4.1}
Let $\xi$ be as in \eqref{eq:2.1}.
For any fixed $B\in(0,C_0)$ there exist two positive constants $\bar{\e}$ and $\bar{\zeta}$ such that
for any $k\ge 1$ and any $\nn\in\ZZZ_{*}^{d}$ one has 
\begin{equation} \nonumber
\left| u^{(k)}_{\nn} \right| \le B_{0} B_1^{k} {\rm e}^{-\xi |\nn|/2} ,
\end{equation}
where $B_{0}$ and $B_{1}$ are two positive constants proportional to $B$,
provided  $|\e|<\bar{\e}$ and $|\zeta|<\bar{\zeta}$.
\end{lemma}

\prova
One bounds \eqref{eq:3.2} as
\begin{eqnarray} \nonumber
\left| \Val_i(\theta) \right| 
& \!\!\! \le \!\!\! &
m^{2k-1} \Biggl( \prod_{v\in V(\theta)} | \e |\Delta \rho^{-p_v} \Biggr) 
\Biggl( \prod_{v \in E_1(\theta)} \| D^{-1}(\oo \cdot \nn_v, \e) \| | \e \, f_{\nn_{v}} | \Biggr) \times
\nonumber \\ & & \qquad \qquad \times 
\Biggl( \prod_{v\in E_{0}(\theta)} |\zeta| \Biggr)
\Biggl( \prod_{v\in V(\theta)} \frac{2}{\ka_{1}^{2}|b_{1}\e|} \Biggr) \nonumber \\
& \!\!\! \le \!\!\! &
m^{2k-1}|\zeta|^{|E_{0}(\theta)|} \Delta^{|V(\theta)|} \rho^{-(|N(\theta)|-1)}
\Phi^{|E_{1}(\theta)|} \times \nonumber \\
& & \qquad \qquad \times 
\Big(\frac{2}{\ka_{1}^{2}|b_{1}|}\Big)^{|V(\theta)|}  
\Biggl( \prod_{v \in E_{1}(\theta)} 
\| D^{-1}(\oo \cdot \nn_v, \e) \| |\e| \, {\rm e}^{-\xi |\nn_{v}|}  \Biggr) , \nonumber
\end{eqnarray}
where we have bounded $f_{\nn_{v}}$ according to \eqref{eq:2.1} and $G_{p_{v}}$ according to \eqref{eq:2.2},
and we have used the bound  $ \| D^{-1}(\oo \cdot \nn_v, \e) \| \le 2/\ka_{1}^{2}|b_{1}\e|$
for the propagators of the lines exiting the nodes $v \in V(\theta)$.
For each end node $v \in E_{1}(\theta)$ we extract a factor ${\rm e}^{-3\xi|\nn_{v}|/4}$,
so that, if we define $C_0$ as in \eqref{eq:3.6}, we obtain
\begin{equation} \label{eq:4.1}
\left| \Val(\theta) \right| \le \rho \, C_{0}^k |\zeta|^{|E_{0}(\theta)|}
\Biggl( \prod_{v\in E_{1}(\theta)}{\rm e}^{-3\xi |\nn_{v}|/4} \Biggr)
\Biggl( \prod_{v \in E_{1}(\theta)} \| D^{-1}(\oo \cdot \nn, \e) \|  |\e| {\rm e}^{-\xi |\nn_{v}|/4}  \Biggr) .
\end{equation}

By the definition \eqref{eq:3.6a}, for any given $N \in\NNN$ we have
$|\oo\cdot\nn| \ge s_{N}$ for all $\nn\in\ZZZ^d_*$
such that $|\nn|\le N$. Set
\begin{equation} \nonumber
\de=\de(N):={\rm e}^{-\xi N/4} .
\end{equation}
Let $B$ be such that $0<B<C_0$. We first fix $N$ such that $2\delta C_{0}^{2} \le 
B^2 \ka_{1}^{2}|b_{1}|$, then we fix
$\bar{\e}$ and $\bar{\zeta}$ by requiring that $C_0^2 \bar{\e}/\ka_1^2 r_N\le B^2$
and $C_0^2\bar{\zeta}<B^2$.

By \eqref{eq:3.5}, in \eqref{eq:4.1}, for all $v\in E_{1}(\theta)$,
we can bound  $\| D^{-1}(\oo\cdot\nn_{v},\e)|| \le 1/\ka_1^2 r_N$ if $|\nn_{v}|\le N$ and
$\| D^{-1}(\oo\cdot\nn_{v},\e)|| \le 2/\ka_1^2|b_{1}\e|$ if $|\nn_{v}|> N$.  Thus,
for all $v \in E_{1}(\theta)$, one has
\begin{equation} \label{eq:4.2}
 |\e| {\rm e}^{-\xi |\nn_{v}|/4} \| D^{-1}(\oo\cdot\nn_{v},\e) \| \le \frac{1}{\ka_1^2}
\max\left\{ \frac{2\delta}{|b_1|} , \frac{|\e|}{r_N} \right\} ,
\end{equation}
so that in \eqref{eq:4.1}, when $2\delta/|b_{1}| \ge |\e|/r_N$, we can bound
\begin{equation} \nonumber
C_{0}^k |\zeta|^{|E_{0}(\theta)|} \Biggl( \prod_{v \in E_{1}(\theta)}  
 |\e| \, {\rm e}^{-\xi |\nn_{v}|/4}\|D^{-1}(\oo\cdot\nn_{v},\e)\| \Biggr) \le
C_{0}^{k} \bar{\zeta}^{|E_{0}(\theta)|} \left( \frac{2\delta}{|b_{1}|\ka_{1}^{2}} 
\right)^{|E_{1}(\theta)| } \le \frac{B}{C_0}B^{k} ,
\end{equation}
while, when $2\delta/|b_{1}| < |\e|/r_N$, we can bound
\begin{equation} \nonumber 
C_{0}^k |\zeta|^{|E_{0}(\theta)|}
\Biggl( \prod_{v \in E(\theta)}   |\e| \, {\rm e}^{-\xi |\nn_{v}|/4}
\|D^{-1}(\oo\cdot\nn_{v},\e)\|  \Biggr)
\le C_{0}^{k} \bar{\zeta}^{|E_{0}(\theta)|}
\left( \frac{|\e|}{\ka_{1}^{2}r_N} \right)^{|E_{1}(\theta)| } \le \frac{B}{C_0} B^{k} ,
\end{equation}
where we have used twice \eqref{eq:3.4} with $k(\theta)=k$.

Summarising, for all $k\in\NNN$ and all $\nn\in\ZZZ^{d}_{*}$ we have obtained
\begin{equation} \nonumber
\left| \Val(\theta) \right| \le B_{0} 
B^{k} \Biggl( \prod_{v\in E(\theta)}{\rm e}^{-3\xi |\nn_{v}|/4} \Biggr) ,
\qquad B_{0} := \rho \frac{B}{C_{0}} .
\end{equation}
The latter bound can be used to provide an estimate for the coefficients $u^{(k)}_{\nn}$. 
According to \eqref{eq:3.3}, we have to sum over all trees in $\TT_{k,\nn}$: 
the sum over  the Fourier labels $\{\nn_{v}\}_{v\in E_{1}(\theta)}$ is performed by using
the factors ${\rm e}^{-3\xi |\nn_{v}|/4}$ associated with the end nodes in $E_{1}(\theta)$,
and gives a bound $C_{1}^{|E_{1}(\theta)|}{\rm e}^{-\xi |\nn|/2}$,
for some positive constant $C_{1}$. 
The sum over the other labels produces a factor $C_{2}^{|N(\theta)|}$,
with $C_{2}$ another suitable positive constant. Set $B_1:=BC_{1}C_{2}$;
then the assertion follows.
\EP


\subsection{Analyticity in $\boldsymbol{\pps}$}
\label{sec:3.5}

The bounds on the coefficients $u^{(k)}_{\nu}$ given in Lemma \ref{lem:4.1} assure
that the series \eqref{eq:2.8}, with $\mu=1$, converges to a well defined function $u(\om,\e,\zeta)$,
which depends analytically on $\psi=\om t$. This is implied by the following result.

\begin{lemma}  
\label{lem:4.3}
For any $\oo\in\RRR^{d}$ there exist two positive constants $\bar{\e}$ and $\bar{\zeta}$
such that the series \eqref{eq:2.8}
converges to a function $u(\pps,\e,\zeta):=u(\pps,\e,\zeta,1)$, with the following properties:
\begin{enumerate}
\itemsep0em
\item $u(\pps,\e,\zeta)$ is analytic in $\pps$ in a strip $\Sigma_{\x'}$, with $\xi'<\xi/2$,
\item $u(\oo t,\e,\zeta)$ solves \eqref{eq:2.6a},
\end{enumerate}
for $\mu=1$, $|\e|<\bar{\e}$ and $|\zeta|<\bar{\zeta}$.
\end{lemma}

\prova
In Lemma \ref{lem:4.1} we can fix $B$ so that $B_1<1$. 
This implies the convergence of the series \eqref{eq:2.8}
provided $\mu$ is taken to satisfy $B_1\mu<1$, which allows $\mu=1$.
The function (\ref{eq:2.8}) solves (\ref{eq:2.7}) order by order by construction, and,
since the series converges uniformly, then it is also a solution \textit{tout court}
to \eqref{eq:2.7} with $\mu=1$ and hence of (\ref{eq:2.6a}). 

The bound on the Fourier coefficients given by Lemma \ref{lem:4.1} assures that
the solution is analytic in $\pps\in\Sigma_{\xi'}$ for any $\xi'<\xi/2$.
\EP

\subsection{Continuity in $\boldsymbol{\e}$}
\label{sec:3.6}

The bounds of Lemma \ref{lem:4.1} are not enough to prove that the solution $u(\om t,\e,\zeta)$
vanishes in the limit $\e\to0$. Indeed, the constants $B$ and $B_0$ are independent of $\e$.
However, a slight adaptation of the proof of Lemma \ref{lem:4.1} leads to the following result.

\begin{lemma} \label{lem:5.1}
Let $\bar{\e}$ and $\bar{\zeta}$ be as in Lemma \ref{lem:4.3}. For any $|\zeta|<\bar{\zeta}$,
the function $u(\psi,\e,\zeta)$  is continuous in $\e\in[0,\bar{\e})$
and it tends to $0$ as $\e\to0$.
\end{lemma}

\prova
Let $\zeta$ be such that $|\zeta|<\bar{\zeta}$. Continuity in $\e$ is obvious for $\e>0$.
To the study the continuity at $\e=0$, set
\begin{equation} \nonumber
F(\e,\zeta) = \|u(\cdot,\e,\zeta)\|_{\io} := \sup\{ u(\pps,\e,\zeta) : \pps \in \Sigma_{\xi'} \} ,
\end{equation}
with $\xi'$ as in Lemma \ref{lem:4.3}. Since $F(0,\zeta)=0$,
we have only to prove that $F(\e,\zeta)  \to 0$ as $\e\to 0$,
that is that for all $\eta>0$ there exists $\e_1>0$ such that $0<\e<\e_1$ implies $|F(\e,\zeta)|<\eta$.
	
Let $\bar{\e}$ be as in Lemma \ref{lem:4.3}. By writing $u(\oo t,\e,\zeta)$ as in \eqref{eq:2.8} with $\mu=1$,
we obtain
%
%
%
\begin{equation} \nonumber
F(\e,\zeta) \le \sum_{k=1}^{\io} \sum_{\nn\in\ZZZ^{d}} \bigl| u^{(k)}_{\nn} \bigr| \, {\rm e}^{\x'|\nn|} ,
\end{equation}
where $\xi'<\xi/2$. By reasoning as in the proof of Lemma \ref{lem:4.1}, we can bound
\begin{eqnarray} \nonumber
\sum_{\nn\in\ZZZ^{d}} \bigl| u^{(k)}_{\nn} \bigr| \, {\rm e}^{\x'|\nn|}
& \!\!\!\! \le \!\!\!\! &
\sum_{\nn\in\ZZZ^{d}} \sum_{\theta \in \TT_{k,\nn}} \bigl| \Val(\theta) \bigr| \, {\rm e}^{\x'|\nn|} \nonumber \\
& \!\!\!\! \le \!\!\!\! &
\rho \, C_{0}^{k} \sum_{\nu \in \ZZZ^{d}} {\rm e}^{\frac{\xi}{2} |\nu|}
\sum_{\theta \in \TT_{k,\nn}} \!\! |\zeta|^{|E_{0}(\theta)|} \times
\nonumber \\
& & \qquad \qquad \times
\Biggl(\prod_{v \in E_{1}(\theta)} e^{-\frac{3}{4}\xi |\nu_{v}|} \Biggr)
\Biggl(\prod_{v \in E_{1}(\theta)} \|D^{-1}(\omega \cdot \nu_{v},\e)\| \, |\e| \,
{\rm e}^{-\frac{\xi}{4} |\nu_{v}|} \Biggr) \nonumber \\
& \!\!\!\! \le \!\!\!\! &
\rho \, C_{0}^{k} \sum_{\nu \in \ZZZ^{d}} e^{\frac{\xi}{2} |\nu|} \sum_{\theta \in \TT_{k,\nn}} \!\! |\zeta|^{|E_{0}(\theta)|}
\Biggl( \frac{1}{\ka_{1}^{2}} \max\bigg\{\frac{2\delta}{|b_{1}|} , \frac{|\e|}{r_N}\bigg\}
\Biggr)^{|E_{1}(\theta)|-1} \!\!\!\!\!\!\!\!\!\! \!\!\!\!\!\!\!\!\!\!  \times \nonumber \\ 
& & \qquad \qquad \times \,
\Biggl(\prod_{v \in E_{1}(\theta)\setminus\{\bar v\}} {\rm e}^{-\frac{3}{4}\xi |\nu_{v}|} \Biggr)
\sum_{\nn_{\bar v}\in\ZZZ^{d}_{*}} |\e| \, {\rm e}^{-\frac{\xi}{4}|\nu_{\bar v}|} \| D^{-1}(\omega \cdot \nu_{\bar v},\e) \| \nonumber \\
& \!\!\!\! \le \!\!\!\! & 
\rho \, C_{0} B^{-1} B_1^{k} \!\! \sum_{\nn \in \ZZZ^{d}_*} |\e|\,{\rm e}^{-\frac{\xi}{4} |\nn|}
\| D^{-1}(\oo\cdot\nn,\e) \| , \nonumber
\end{eqnarray}
with $B_1$ as in Lemma \ref{lem:4.1}. We have used that any tree $\theta$ contains
at least one end node $\bar v$ with $\nu_{\bar v}\neq 0$, since $\nn\neq\vzero$.
In particular, this implies $|E_{1}(\theta)|\ge 1$. 
Therefore we may use the bound \eqref{eq:4.2} for all the end nodes $v\in E_{1}(\theta)$ except for one.
As a consequence, we obtain, for any $N\in\NNN$,
\begin{equation}  \label{eq:5.1}
F(\e,\zeta) \le \frac{\rho \, C_0 \,B^{-1}}{1-B_1} \sum_{\substack{\nn\in\ZZZ^d \\ |\nn| \le N}}
|\e|\,{\rm e}^{-\xi |\nn|/4} \| D^{-1}(\oo\cdot\nn,\e) \| + D_0 {\rm e}^{-\xi N/8} ,
\end{equation}
where we have set
\begin{equation} \nonumber
D_0 := \frac{\rho \, C_0 \, B^{-1}}{1-B_1} \sum_{\nn \in \ZZZ^{d}} {\rm e}^{-\xi|\nn|/8} .
\end{equation}
Fix $\eta>0$ and take $N_1$ such that $D_0 {\rm e}^{-\xi N_1/8} < \eta/2$. If we define
$r_{N}$ as in \eqref{eq:2.6a} and choose $\e_1$ small enough, we bound
\begin{equation} \nonumber
\frac{\rho \, C_0 \, B^{-1}}{1-B_1} \sum_{\substack{\nn\in\ZZZ^d \\ |\nn| \le N_1}}|\e|\,{\rm e}^{-\xi |\nn|/4}
\| D^{-1}(\oo\cdot\nn,\e) \| \le \frac{|\e|\, D_{0} }{r_{N_1}} < \frac{\eta}{2},
\end{equation}
provided one has $0<\e<\e_1$. This implies $|F(\e,\zeta)| < \eta$. Hence the thesis follows.
\EP


\zerarcounters 
\section{The bifurcation equation}
\label{sec:4}

In this section we show that, given the solution $u(\oo t,\e,\zeta)$ to the range equation
\eqref{eq:2.6a}, it is possible to fix $\zeta$ as a function of $\e$ in such a way that
the bifurcation equation \eqref{eq:2.6b} is solved as well.

\begin{lemma} \label{lem:5.2}
Let $\bar{\e}$ and $\bar{\zeta}$ be as in Lemma \ref{lem:4.1} and let $u=u(\oo t,\e,\zeta)$
be the function in Lemma \ref{lem:4.3}. There exist neighbourhoods $U \subset(-\bar{\e},\bar{\e})$
and $V \subset B_{\bar\zeta}(0)$ and a function $\zeta : U \to V$ with the following properties:
\begin{enumerate}
\itemsep0em
\item the function $\e \mapsto \zeta(\e)$ is continuous in $U$,
\item for all $\e\in U$ the equation \eqref{eq:2.6b} is satisfied for $\zeta=\zeta(\e)$,
\item $\zeta(\e)$ is the only solution to \eqref{eq:2.6b} in $V$.
\end{enumerate}
\end{lemma}

\prova
Write \eqref{eq:2.6b} as
\begin{equation} \nonumber
\calH(\zeta,\e) := A \, \zeta + \left[ G(c+\zeta+u) \right]_{0} = 0 .
\end{equation}
By construction, the function $u(\pps,\e,\zeta)$ depends analytically on $\zeta$
in a neighbourhood of the origin. Therefore,
$\calH(\zeta,\e)$ is analytic in $\zeta$ and, by Lemma \ref{lem:5.1},  is continuous in $\e$. One has
$\calH(0,0) = 0$ and
\begin{equation} \nonumber
\frac{\partial \calH_{i}}{\partial \zeta_{j}}(0,0) = \frac{\partial}{\partial \zeta_{j}}(A \zeta + G(c + \zeta + u ))_{i}(0,0)
= A_{i,j} , \qquad i,j=1,\dots, m.
\end{equation}
Since $\det A \ne 0$ by Hypothesis \ref{hyp:1}, we can apply the implicit function theorem, in the version of
Loomis and Sternberg \cite{LS},  so as to conclude that there exist neighbourhoods 
$U\subset(-\bar{\e},\bar{\e})$ and $V\subset B_{\bar{\zeta}}(0) $
such that for all $\e\in U$ one can find a unique value $\zeta(\e) \in V$,
depending continuously on $\e$, such that $\calH(\zeta(\e),\e)=0$.
\EP

\begin{lemma} \label{lem:5.3}
Let the functions $u(\pps,\e,\zeta)$ and $\zeta(\e)$ be as defined in Lemma \ref{lem:4.3}
and in Lemma \ref{lem:5.2}, respectively.
There exists $\e_{0}>0$ such that for all $\e\in(0,\e_0)$ the function 
$x(t,\e)=c+\zeta(\e)+u(\oo t,\e,\zeta(\e))$ solves \eqref{eq:2.3}.
Moreover $x(t,\e) \to c$ as $\e\to 0$.
\end{lemma}

\prova
The result follows immediately from Lemma \ref{lem:4.3} and Lemma \ref{lem:5.2}.
\EP


\zerarcounters 
\section{Extension of the results to more general nonlinearities}
\label{sec:6}

In this section, we show how to extend the results of the previous section,
so as to prove Theorem \ref{thm:2}.

\subsection{Recursive equations}

Let the function $h(x,\pps)$ in \eqref{eq:1.4} be analytic on the domain
$B_{\rho_0}(c) \times \Sigma_{\xi_0}$, with $B_{\rho_0}(c)$ and  $\Sigma_{\xi_0}$ as in Section \ref{sec:2}.
We define $H_{p}(c,\psi)$, the tensor of rank $p+1$ with components
\begin{equation} \nonumber
[H_{p}(c,\psi)]_{i, i_{1},\dots, i_{p}} = \left. \frac{1}{p!}
\frac{\partial^p h_{i}}{\partial x_{i_{1}},\dots, \partial x_{i_{p}}}(x,\psi) \right|_{\xx=\cc} ,
\qquad i=1,\ldots, m, 
\end{equation}
where $H_{0}(\psi,c)=h(\psi,c)$ and $i_{1}, \dots, i_{p}= 1, \dots, m$ for $p\ge 1$, and write
\begin{equation} \nonumber
h_{i}(x,\psi) = \sum_{p=0}^{\infty} [H_{p}(c,\psi)(x-c)^{p}]_{i} = 
\sum_{p=0}^{\infty}\sum_{i_{1},\dots,i_{p}=1}^{m} [H_{p}(c,\psi)]_{i,i_{1},\dots,i_{p}} 
\prod_{j=1}^{p}(x_{i_{j}}-c_{i_{j}}) .
\end{equation}
Then we consider the Fourier expansions
\begin{equation} \nonumber
h(x,\pps) = \sum_{\nn\in\ZZZ^{d}} {\rm e}^{i\nn\cdot\pps} \, h_{\nn}(x) ,
\qquad H_{p}(c,\pps) = \sum_{\nn\in\ZZZ^{d}} {\rm e}^{i\nn\cdot\pps} \, H_{p,\nn}(c) ,
\end{equation}
so that, for any positive $\rho<\rho_0$ and $\xi<\xi_0$,
\begin{equation} 
\label{eq:6.0}
\left| [H_{p,\nu}(c)]_{i,i_{1},\dots,i_{p}} \right| \le \Delta \rho^{-p} {\rm e}^{-\xi|\nn|} ,
\end{equation}
for a suitable positive constant $\Delta$. We rewrite \eqref{eq:1.4} as
\begin{eqnarray}
\e M \ddot x + \Gamma \dot{x}   
& \!\!\! + \!\!\! & 
\e \, A \left( x - c \right) + \e \sum_{\nn\in\ZZZ^{d}_{*}} {\rm e}^{i\nn\cdot\oo t} h_{\nn}(c) \nonumber \\
& \!\!\! + \!\!\! &  
\e \sum_{\nn\in\ZZZ^{d}_{*}} H_{1,\nn}(c) \, {\rm e}^{i\nn\cdot\oo t} \left( x - c \right) +
\e \sum_{p=2}^{\io} \sum_{\nn\in\ZZZ^{d}} H_{p,\nn}(c) \, {\rm e}^{i\nn\cdot\oo t} 
\left( x - c \right)^{p} = 0,\nonumber 
\end{eqnarray}
where we have used that $h_{\vzero}(c)=0$ and the $m \times m$ matrix  $A:= H_{1,0}(c)$, with entries 
$$A_{i,j}=\left. \frac{\partial h_{i,0}( x)}{\partial x_{j}} \right |_{\xx=\cc} = 
\left. \frac{\partial^{2} \tilde{V}_{0} (x)}{\partial x_{i}\partial x_{j}} \right |_{\xx=\cc} , \qquad i,j=1, \dots, m , $$
is positive definite by Hypothesis \ref{hyp:2}.

We look for a solution of the form \eqref{eq:2.4}. If we pass to Fourier space, setting
\begin{equation} \nonumber
\al_1 = \al_{1}(c,\pps) := \sum_{\nn\in\ZZZ^{d}_{*}} H_{1,\nu}(c) \, {\rm e}^{i\nn\cdot\pps} , \qquad
\al_p = \al_{p}(c,\pps) := \sum_{\nn\in\ZZZ^{d}} H_{p,\nu}(c) \, {\rm e}^{i\nn\cdot\pps} ,
\end{equation}
and defining  $D(\e,s)$ as (compare with \eqref{eq:2.5})
\begin{equation} \label{eq:6.1}
D(\e,s):= -\e s^{2} M + \ii s \Gamma  + \e A ,
\end{equation}
we obtain the equations
\begin{subequations} \label{eq:6.2}
	\begin{align}
	D(\e,\oo\cdot\nn) \, u_{\nn} & = - \e h_{\nn}(c) - 
	\e \left[ \al_{1} \left( x-c \right) \right]_{\nn} - \e \sum_{p=2}^{\io}
	\left[ \al_{p} \left( x-c \right)^{p} \right]_{\nn} , \qquad \nn \neq \vzero ,
	\label{eq:6.2a} \\
	\e \, A \, \zeta  & = - \left[ \al_{1} \left( x-c \right) \right]_{\vzero} - \e \sum_{p=2}^{\io}
	\left[ \al_{p} \left( x-c \right)^{p} \right]_{\vzero} .
	\label{eq:6.2b}
	\end{align}
\end{subequations}

By following the same strategy as in Section \ref{sec:3}, we shall study first
\eqref{eq:6.2a} and look for a solution depending on the parameter $\zeta$.
Thereafter we shall fix $\zeta$ by requiring that \eqref{eq:6.2b} is solved as well.

Instead of \eqref{eq:6.2a} we consider the equation
\begin{equation} 
\label{eq:6.3}
D(\e,\oo\cdot\nn) \, u_{\nn} = - \mu \, \e h_{\nn}(c) -  \mu \, \e
\left[ \al_{1} \left( x-c \right) \right]_{\nn} - \mu \, \e\sum_{p=2}^{\io}
\left[ \al_{p} \left( x-c \right)^{p} \right]_{\nn} , \qquad \nn \neq \vzero ,
\end{equation}
and look for a quasi-periodic solution to (\ref{eq:6.3}) in the form of a power series in $\mu$,
\begin{equation} 
\label{eq:6.4}
x(t,\e,\mu) = c + \zeta + u(\oo t, \e,\zeta,\mu) , \qquad
u(\oo t, \e,\zeta,\mu) = \sum_{k=1}^{\io} \sum_{\nn\in\ZZZ^{d}_{*}} 
\mu^{k} {\rm e}^{\ii\nn\cdot\pps} u^{(k)}_{\nn} .
\end{equation}
We find the recursive equations
\begin{subequations} \label{eq:6.5}
	\begin{align}
	D(\e, \oo\cdot\nn) \, u^{(1)}_{\nn} & = - \e \, h_{\nn}(c) 
	\label{eq:6.5a} \\
	D(\e, \oo\cdot\nn) \, u^{(k)}_{\nn} & = - \e \, \sum_{\nn_0 \in\ZZZ^{d}_{*}}
	H_{1,\nn_{0}}(c) \, u^{(k-1)}_{\nn-\nn_0} 
	\nonumber \\
	&  - \e \, \sum_{p=2}^{\io} \!\!\!\!\!\!
	\sum_{\substack{ k_{1},\ldots,k_{p} \ge 1 \\ k_{1}+\ldots+k_{p}=k-1}}
	\sum_{\substack{\nn_{0},\nn_{1},\ldots,\nn_{p} \in\ZZZ^{d} \\ \nn_{0}+\nn_{1}+\ldots+\nn_{p}=\nn}} \!\!\!\!\!\!
	H_{p, \nn_{0}}(c) \, u^{(k_{1})}_{\nn_{1}} \ldots u^{(k_{p})}_{\nn_{p}} , \qquad k \ge 2 ,
	\label{eq:6.5b}
	\end{align}
\end{subequations}
where $\nn\neq\vzero$ and we have set once more $u^{(1)}_{\vzero}=\zeta$ and $u^{(k)}_{\vzero}=0$ $\forall k\ge 2$.

\subsection{Trees and chains}
As in Section \ref{sec:3}, we can represent the coefficients $u^{(k)}_{\nn}$ with tree diagrams; there are just few differences with the previous case. 
Define the sets $N(\theta)$, $E(\theta)$, $E_{0}(\theta)$, $E_{1}(\theta)$, $V(\theta)$ and $L(\theta)$ as previously
and call $k(\theta):=|N(\theta)|$ the order of $\theta$.
Now we split $V(\theta)=V_{1}(\theta) \amalg V_{2}(\theta)$, with
$V_{1}(\theta)=\{ v \in V(\theta) : p_{v} = 1 \}$ and $V_{2}(\theta)=\{ v \in V(\theta) : p_{v} \geq 2 \}$.
By contrast with Section \ref{sec:3}, now in general $V_{1}(\theta) \neq \emptyset$.

We associate with each node $v\in N(\theta)$ a \textit{mode} label $\nn_{v}\in\ZZZ^{d}$
and with each line $\ell\in L(\theta)$ a \textit{momentum} $\nn_{\ell} \in \ZZZ^{d}$ defined by
\begin{equation} \nonumber  
\nn_{\ell}=\sum_{\substack{w\in N(\theta) \\ w \prec \ell}} \nn_{w} .
\end{equation}
Finally we impose the constraints that
\begin{itemize}
	\itemsep0em
	\item $\nn_{v} \neq \vzero$ $\forall v\in V_{1}(\theta)$,
	\item $\nn_{\ell} \neq \vzero$ for any line $\ell$ exiting a node in $V(\theta)$.
\end{itemize}

We associate with each node $v\in N(\theta)$  a \textit{node factor}
\begin{equation} \label{6.6}
F_{v} := \begin{cases}
- \e \, H_{p_{v},\nn_{v}}(c) , & v \in V(\theta) , \\
-\e \, h_{\nn_{v}}(c) , & v \in E_{1}(\theta) , \\
\zeta , & v \in E_{0}(\theta) ,
\end{cases}
\end{equation}
and with each line $\ell\in L(\theta)$ a \textit{propagator}
\begin{equation} \label{eq:6.7}
\calG_{\ell} := \begin{cases}
D^{-1}(\e, \oo\cdot\nn_{\ell}) , & \nn_{\ell} \neq \vzero , \\
\uno , & \nn_{\ell}=\vzero .
\end{cases}
\end{equation}
We define the value of the tree $\theta$ as in \eqref{eq:3.2} and write
$u^{(k)}_{\nn}$ as in \eqref{eq:3.3}, where $\TT_{k,\nn}$ denotes the set of non-equivalent trees
of order $k$ and momentum $\nn$ associated with the root line, constructed according to the
new rules given above. Then, by construction, the coefficients $u^{(k)}_{\nn}$ solve formally \eqref{eq:6.5}.

The main difference with respect to the trees considered in Section \ref{sec:3}
is that, now, trees may contain `chains'.
In a tree $\theta$ we define a \emph{chain} $\CCCC$ as a subset of $\theta$ formed
by a connected set of nodes $v \in V(\theta)$ with $p_{v}=1$ and by the lines exiting them.
Therefore, if $V(\CCCC)$ and $L(\CCCC)$ denote the set of nodes and the set of lines of $\CCCC$ and 
$V(\CCCC)=\{v_1,v_2,\ldots,v_p\}$, with $v_1 \succ v_2 \succ \ldots \succ v_p$,
then $L(\CCCC)=\{\ell_{1},\ell_{2},\ldots,\ell_{p}\}$, where
$\ell_{i}$ is the line exiting $v_{i}$, for $i=1,\ldots,p$, and $\ell_{i}$ enters the node $v_{i-1}$ for $i=2,\ldots,p$.
We call $p=|V(\CCCC)|=|L(\CCCC)|$ the \emph{length} of the chain $\CCCC$. Finally we define
the \emph{value} of the chain $\CCCC$ as the matrix
\begin{equation} \nonumber
\Val(\CCCC) := \Biggl( \prod_{v\in V(\CCCC)} F_{v} \Biggr) 
\Biggl( \prod_{\ell\in L(\CCCC)} \calG_{\ell} \Biggr) ,
\end{equation}
with component indices $i$ and $j$ equal to the first index of the propagator of the line $\ell_1$
and to the first index of the propagator of the line entering $v_p$, respectively,
and denote by $\gotC(\theta)$ the set of all maximal chains contained in $\theta$.

\subsection{Bounds}

We can see that the bounds on the propagators are the same as in Section \ref{sec:3.03}. 
Indeed let $K$ and $\ka_{1} \le \dots \le \ka_{m}$ be as in the aforementioned section, so that
the inverse of $D=D(\e,s)$ can be written as
\begin{equation} \nonumber
D^{-1}=K^{-1 } \left( \ii s \uno + K^{-1} (-\e s^2 M + \e A) K^{-1} \right)^{-1} K^{-1}.
\end{equation}
By defining 
\begin{equation} \nonumber
S := \e \, K^{-1} A K^{-1} - \e s^{2} K^{-1}MK^{-1} \quad \text{and} \quad R:= S + \ii s \uno ,
\end{equation}
we obtain $D^{-1}= K^{-1} R^{-1} K^{-1}$.
Again, let $b_{1} \le \dots \le b_{m}$ be the eigenvalues of $K^{-1}AK^{-1}$:
the eigenvalues of $S$ are $b_{k}\e + O(\e s^{2})$, $k=1,\dots, m$, and hence
the eigenvalues of $R$ are
$b_{k} \e + \ii s + O(\e s^{2})$, $k=1,\dots, m$. 
By repeating the same argument as in Section \ref{sec:3.03}, we find, for $\e$ small enough,
\begin{equation}  \label{eq:6.8}
\| D^{-1} \| \le \frac{1}{\ka^{2}_{1}} \min  \left\{ \frac{2}{|b_{1}\e|} ,
\max \left\{ \frac{1}{\al} , \frac{1}{|s|} \right\} \right\}, 
\end{equation} 
that is the bound that we had in \eqref{eq:3.5}. 

We define $s_N$ and $r_{N}$ as in \eqref{eq:3.6a}, and fix $N$ and $\de$ as in Section \ref{sec:3.4}.
With respect to \eqref{eq:3.6}, we redefine the costant $C_0$ as
\begin{equation} \label{eq:6.9}
C_{0}:=m^2 \rho^{-1}\max\bigg\{\frac{2\Delta}{\ka_{1}^{2}|b_{1}|},1\bigg\} .
\end{equation}
with $\rho$ and $\Delta$ as in \eqref{eq:6.0}, and set
\begin{equation} \label{eq:6.9b}
\be:=\max\bigg\{\de,\frac{2 |\e b_{1}|}{r_{N}}\bigg\} .
\end{equation}
The following lemma provides us a bound for the values of the chains.

\begin{lemma} \label{lem:6.1}
For any $p \ge 1$ and for any chain $\CCCC$ of length $p$, one has
\begin{equation} \nonumber
\| \Val(\CCCC) \| \le C_{0}^{p} \be^{(p-1)/2} \prod_{v\in V(\CCCC)} {\rm e}^{-3\xi|\nn_{v}|/4} .
\end{equation}
\end{lemma}

\prova
We proceed by induction on $p$. 
If $p=1$ then $\CCCC$ contains only one node $v$, so that
\begin{equation} \nonumber
\| \Val(\CCCC) \| \le m |\e| \Delta 
\rho^{-1} {\rm e}^{-\xi|\nn_{v}|} \| D^{-1}(\e,\oo\cdot\nn_{v}) \| \le C_{0} {\rm e}^{-3\xi|\nn_{v}|/4} ,
\end{equation}
where we have bounded $\|D^{-1}(\e,\oo\cdot\nn_{v})\| \le 2/\ka_1^2|b_1\e|$, according to \eqref{eq:6.8}.

Let  us denote by $v_1 \succ v_2 \succ v_3 \succ \ldots \succ v_{p}$ the nodes
in $\CCCC$ and by $\ell_1,\ell_2,\ell_3,\ldots,\ell_p$ the lines exiting such nodes.
If $p\ge 2$, the nodes $\{v_2,v_3,\ldots,v_p\}$ and the lines $\{\ell_2,\ell_3,\ldots,\ell_p\}$
form a chain $\CCCC'$ (contained in $\CCCC$) of length $p-1$; moreover, if $p\ge 3$, the nodes $\{v_3,\ldots,v_p\}$
and the lines $\{\ell_3,\ldots,\ell_p\}$ form a chain $\CCCC''$ (contained in $\CCCC'$) of length $p-2$.

Assume the bound to hold up to $p-1$.
If $|\oo\cdot\nn_{\ell_1}| \ge s_{N}/2$,  one has
\begin{eqnarray} 
\| \Val(\CCCC) \| & \!\!\!\! \le \!\!\!\! & m
|\e| \Delta \rho^{-1} {\rm e}^{-\xi|\nn_{v_1}|} \| D^{-1}(\e,\oo\cdot\nn_{v}) \| \| \Val(\CCCC') \| \nonumber \\
& \!\!\!\! \le \!\!\!\! &
C_{0} \be \left( C_{0}^{p-1} \be^{(p-2)/2} \right) 
\prod_{v\in V(\CCCC)} {\rm e}^{-3\xi|\nn_{v}|/4} , \nonumber
\end{eqnarray}
which yields the bound for $p$.
If $|\oo\cdot\nn_{\ell_1}| < s_{N}/2$ and $p\ge 3$ we distinguish between two cases.
If $|\oo\cdot\nn_{\ell_2}| \ge s_{N}/2$, we can bound
\begin{equation} \nonumber
\begin{split}
\| \Val(\CCCC) \| & \le m^{2} |\e| \Delta \rho^{-1} {\rm e}^{-\xi|\nn_{v_1}|} \| D^{-1}(\e,\oo\cdot\nn_{\ell_1}) \| 
|\e| \Delta \rho^{-1} {\rm e}^{-\xi|\nn_{v_2}|} \| D^{-1}(\e,\oo\cdot\nn_{\ell_2}) \| \| \Val(\CCCC'') \| \\
&\le C^2_{0} \be \left( C_{0}^{p-2} \be^{(p-3)/2} \right) 
\prod_{v\in V(\CCCC)} {\rm e}^{-3\xi|\nn_{v}|/4} ,\\
\end{split}
\end{equation}
so that the bound follows once more.

If $|\oo\cdot\nn_{\ell_2}| < s_{N}/2$, then
\begin{eqnarray}
|\oo\cdot\nn_{v_1}| & \!\!\! = \!\!\! & |\oo\cdot\nn_{v_1}+\oo\cdot\nn_{\ell_2} - \oo\cdot\nn_{\ell_2}|
\le |\oo\cdot\nn_{v_1}+\oo\cdot\nn_{\ell_2} | + |\oo\cdot\nn_{\ell_2}| \nonumber \\
& \!\!\! = \!\!\! & |\oo\cdot\nn_{\ell_1} | + |\oo\cdot\nn_{\ell_2}| < s_{N} , \nonumber
\end{eqnarray}
so that, since $\nn_{v_1} \neq \vzero$ and hence $\oo\cdot\nn_{v_{1}} \neq 0$,
we conclude that $|\nn_{v_1}| > N$, which allows us to bound
${\rm e}^{-\xi|\nn_{v_1}|} \le \de {\rm e}^{-3\xi|\nn_{v_1}|/4}$. Therefore we obtain
\begin{equation} \nonumber
\| \Val(\CCCC) \| \le 
C_{0}^{2} \de \| \Val(\CCCC'') \| \le C^2_{0} \de \left( C_{0}^{p-2} \be^{(p-3)/2} \right) 
\prod_{v\in V(\CCCC)} {\rm e}^{-3\xi|\nn_{v}|/4} ,
\end{equation}
which gives the bound for $p$ in this case too.

Finally if $p=2$ we can reason as in the case $p\ge 3$, the only difference being that
the matrix $\Val(\CCCC'')$ has to replaced with $\uno$ -- since there is no further
chain $\CCCC''$ for $p=2$. Therefore we obtain
$\|\Val(\CCCC)\| \le C_{0}^2 \be$, which is the desired bound.
\EP

\begin{lemma} \label{lem:6.2}
For any tree $\theta$ one has $4|E(\theta)| + 2|V_{1}(\theta)| - 2 |\gotC(\theta)| \ge k(\theta) + 2$.

\end{lemma}

\prova
See \cite[page 12]{GV}.
\EP

A result analogous to Lemma \ref{lem:4.1} still holds. This can be proved as follows.
For any tree $\theta$ we have
\begin{eqnarray} \nonumber
\left| \Val_{i}(\theta) \right| 
& \!\!\! \le \!\!\! & 
m^{2k-1-2|V_1(\theta)|} 
\Biggl( \prod_{v\in V_{2}(\theta)} 
\frac{2| \e| \, \Delta \rho^{-p_{v}} {\rm e}^{-\xi |\nn_{v}|}}{\ka_{1}^{2}|\e b_{1}|} \Biggr)
\times \nonumber \\
& \!\!\!  \!\!\! & \times \,
\Biggl( \prod_{v \in E_1(\theta)}  | \e| \Delta \, {\rm e}^{-\xi |\nn_{v}|} \|D^{-1}(\oo\cdot\nn_{v},\e)\| \Biggr) 
\Biggl( \prod_{v\in E_{0}(\theta)} |\zeta| \Biggr) 
\Biggl( \prod_{\CCCC \in \gotC(\theta)} \| \Val(\CCCC) \|  \Biggr) \nonumber \\
& \!\!\! \le \!\!\! & 
m^{-1} \rho \, C_{0}^{k} |\zeta|^{|E_{0}(\theta)|} 
\Biggl( \prod_{v \in E_{1}(\theta)} 
|\e| \, {\rm e}^{-\xi |\nn_{v}|/4} \|D^{-1}(\oo\cdot\nn_{v},\e)\| \Biggr) \times \nonumber \\
 & \!\!\!  \!\!\! & \times \,
\Biggl( \prod_{\CCCC \in \gotC(\theta)} \be^{(|V(\CCCC)|-1)/2} \Biggr)
\Biggl( \prod_{v \in V(\theta)} {\rm e}^{-3\xi |\nn_{v}|/4} \Biggr) , \nonumber
\end{eqnarray}
for $i=1, \dots, m$, so that
\begin{equation} \nonumber
\left| \Val(\theta) \right| \le \rho \, C_{0}^{k} \bar{\zeta}^{|E_{0}(\theta)|} 
\left( \frac{1}{\ka_{1}^{2}} \max \left\{ \frac{2\de}{|b_{1}|},\frac{\bar{\e}}{r_{N}} \right\} \right)^{|E_{1}(\theta)|}
\be^{(|V_1(\theta)|-|\gotC(\theta)|)/2} 
\Biggl( \prod_{v \in V(\theta)} {\rm e}^{-3\xi |\nn_{v}|/4} \Biggr) .
\end{equation}
Therefore, for any constant $B\in(0,C_0)$, if we fix first $N$ and hence $\bar{\e}$ and $\bar{\zeta}$
so that for any $|\e|<\bar{\e}$ and any $|\zeta|<\bar{\zeta}$ one has
\begin{equation} \label{eq:6.10}
C_{0}^{4} \max 
\left\{ \bar{\zeta}, \frac{2\de}{\ka^{2}_{1} |b_{1}|}, \frac{\bar{\e}}{\ka^{2}_{1} r_{N}} ,\beta \right\} < B^4 ,
\end{equation}
we obtain
\begin{equation} \nonumber 
\left| \Val(\theta) \right| \le \rho \, C_{0}^{k}
\left( \frac{B}{C_{0} }\right)^{4|E(\theta)| + 2|V_{1}(\theta)| - 2 |\gotC(\theta)|}
\Biggl( \prod_{v \in V(\theta)} {\rm e}^{-3\xi |\nn_{v}|/4} \Biggr) .
\end{equation}
By using Lemma \ref{lem:6.2} we arrive at
%
%
%
%
\begin{equation} \nonumber 
\left| \Val(\theta) \right| \le \rho \, C_{0}^{k}
\left( \frac{B}{C_{0} }\right)^{k+2} \Biggl( \prod_{v \in V(\theta)} {\rm e}^{-3\xi |\nn_{v}|/4} \Biggr) \le B_{0}  B^{k} 
\Biggl( \prod_{v \in V(\theta)} {\rm e}^{-3\xi |\nn_{v}|/4} \Biggr) ,
\qquad B_{0} := \rho \frac{B^2}{C_{0}^2} .
\end{equation}
From here on, we can reason as in Section \ref{sec:3} and we find that the coefficients
$u^{(k)}_{\nn}$ can be bounded as
\begin{equation} \nonumber
\left| u^{(k)}_{\nn} \right| \le B_{0} B_2^{k} {\rm e}^{-\xi |\nn|/2} ,
\end{equation}
for a suitable constant $B_2$ proportional to the constant $B$ in \eqref{eq:6.10}.

\subsection{Conclusion of the proof}

For $\bar{\e}$ and $\bar{\zeta}$ small enough, $B_2$ can be made arbitrarily small;
in particular one obtains $B_2 \le \tilde{B}$, with $\tilde B <1$, so that the series \eqref{eq:6.4}
converges for $\mu=1$. As a consequence,
the function $x(t,\e,\zeta)$ solves the range equation \eqref{eq:6.2a} for any $\zeta$ small enough.

By reasoning as in Section \ref{sec:3.6} and \ref{sec:4}, respectively,
one proves the continuity in $\e$ of the function $x(t,\e,\zeta)$
and fixes the parameter $\zeta$ when dealing with the bifurcation equation \eqref{eq:6.2b}.
This concludes the proof of Theorem \ref{thm:2}.



\end{document}